\documentclass[12pt]{article}

\title{Introduction to work of Hassett-Pirutka-Tschinkel and Schreieder}
\author{Jean-Louis Colliot-Th\'el\`ene}
\date{April  15th,  2018}
\usepackage{color}
\usepackage{amssymb}
\usepackage{amsmath}
\usepackage{amscd}
\usepackage[all]{xy}

\DeclareFontFamily{U}{wncy}{}
\DeclareFontShape{U}{wncy}{m}{n}{%
   <5>wncyr5%
   <6>wncyr6%
   <7>wncyr7%
   <8>wncyr8%
   <9>wncyr9%
   <10>wncyr10%
   <11>wncyr10%
   <12>wncyr6%
   <14>wncyr7%
   <17>wncyr8%
   <20>wncyr10%
   <25>wncyr10}{}
\DeclareMathAlphabet{\cyr}{U}{wncy}{m}{n}

\newtheorem{theo}{Theorem}[section]
\newtheorem{prop}[theo]{Proposition}
\newtheorem{lem}[theo]{Lemma}
\newtheorem{cor}[theo]{Corollary}
\newtheorem{defi}[theo]{Definition}

\newcommand{\bthe}{\begin{theo}}
\newcommand{\ble}{\begin{lem}}
\newcommand{\bpr}{\begin{prop}}
\newcommand{\bco}{\begin{cor}}
\newcommand{\bde}{\begin{defi}}
\newcommand{\ethe}{\end{theo}}
\newcommand{\ele}{\end{lem}}
\newcommand{\epr}{\end{prop}}
\newcommand{\eco}{\end{cor}}
\newcommand{\ede}{\end{defi}}

\def \Gal {{\rm{Gal}}}

\def \Spec {{\rm{Spec}}}
\def\et{{\rm{\acute et}}}

\def \Br {{\rm{Br}}}

\def \to {{\rightarrow}}

\def \A {{\mathbb A}}
\def \Q {{\mathbb Q}}
\def \Z {{\mathbb Z}}

\def \C {{\mathbb C}}

\def \G {{\mathbb G}}
\def \P {{\mathbb P}}

\usepackage[backref=page]{hyperref}
\hypersetup{colorlinks=true,linkcolor=blue,anchorcolor=blue,citecolor=blue}

\begin{document}

\maketitle

Hassett, Pirutka and Tschinkel \cite{HPT} gave the first examples of
families $X \to B$ of smooth, projective,  connected, complex varieties
having some rational fibres
and some other fibres which are not even stably rational.  This used the specialisation
method of Voisin, as extended by   Pirutka and myself.
Under specific circumstances, a simplified version
 of the specialisation method was
produced by Schreieder \cite{Schr1, Schr2}, leading   to a simpler proof   of the HPT example
(no explicit resolution of singularities). In the following note, written on
  the occasion of 
the conference  Quadratic Forms in Chile 2018, held at IMAFI, Universitad de Talca, 8-12 January 2018,  
 I  describe the method in its simplest form. 
 For further developments,  the reader is invited to read \cite{ABBP}, which offers a different look at 
 \cite{HPT} as well as some generalizations, 
   \cite{ABP}, and the papers
  \cite{Schr1,Schr2,Schr3} by Schreieder.
    
    I thank Asher Auel for remarks on the typescript.

\section{Basics on the Brauer group and on the Chow group of zero-cycles}

Grothendieck defined  the Brauer group $\Br(X)$ of a scheme $X$ as
the second \'etale cohomology group $H^2_{\et}(X, \G_{m})$ 
of $X$ with values in the sheaf $\G_{m,X}$ on $X$.
This is a contravariant functor with respect to arbitrary morphisms
of schemes.

If $X=\Spec(K)$ is the spectrum of a field, then
$\Br(X)=\Br(K)$, the more classical cohomological Brauer group
$H^2(\Gal(K_{s}/K),K_{s}^*)$.
Assume $2\in K^*$.  To the quaternion algebra $(a,b)$ one associates
a class $(a,b)\in \Br(K)[2]$. The quaternion algebra is isomorphic to a matrix
algebra $M_{2}(K)$ if and only $(a,b)=0 \in \Br(K)$, if and only if the
diagonal quadratic form $<1,-a,-b>$ has a nontrivial zero over $K$,
if and only if the diagonal quadratic form $<1,-a,-b,ab>$ has a nontrivial zero over $K$.

\begin{prop}
If $X$ is an integral regular scheme and $K(X)$ is its field of rational functions,
then the natural map $\Br(X) \to \Br(K(X))$ is injective.
\end{prop}

\begin{prop}
For $R$ a discrete valuation ring with perfect residue field $\kappa$ and field of fractions $K$,
there is a natural exact sequence
$$0 \to \Br(R) \to \Br(K) \to H^1(\kappa,\Q/\Z) \to 0.$$
The map $\partial_{R}  : \Br(K) \to H^1(\kappa,\Q/\Z)$ is the residue map.
\end{prop}

Suppose $2 \in R^*$. Given $a,b \in K^*$ we may consider the class
$(a,b) \in \Br(K)[2]$ associated to the quaternion algebra $(a,b)$.
Let $v : K^* \to \Z$ be the valuation map.
The quotient $a^{v(b)}/b^{v(a)} \in K^*$ belongs to $R^*$.
One shows :
$$\partial_{R}((a,b)) = (-1)^{v(a).v(b)} cl((a^{v(b)}/b^{v(a)})) \in \kappa^*/\kappa^{*2} =H^1(\kappa,\Z/2) \subset H^1(\kappa,\Q/\Z).$$

\begin{prop}
Let $R \subset S$ be a local inclusion of discrete valuation rings, inducing
an inclusion of fields $K \subset L$ and an inclusion of residue fields $\kappa \subset \lambda$.
Assume $char(\kappa)=0$. Let $e$ be the ramification index.
Then there is a commutative diagram
$$\begin{array}{ccccccccc}
 \Br(K) &\to & H^1(\kappa,\Q/\Z)\\
\downarrow& & \downarrow \\
\Br(L) & \to & H^1(\lambda,\Q/\Z)\\
 \end{array}$$
 where $H^1(\kappa,\Q/\Z) \to H^1(\lambda,\Q/\Z)$ is $e. Res_{\kappa, \lambda}$.
\end{prop}

Let $K$ be a field and $X$ an algebraic variety over $K$, i.e. a $K$-scheme of finite type.
The group  of zero-cycles $Z_{0}(X)$ is the free abelian group on closed points of $X$.
 Given any $K$-morphism $f : Y \to X$ of $K$-varieties,  one defines
 $f_{*} : Z_{0}(Y) \to Z_{0}(X)$ as the map sending a closed point $M\in Y$ with image the closed point $N=f(M) \in X$ to
 the zero cycle $[K(M):K(N)]N \in Z_{0}(X)$.
 
Given a normal, connected curve $C$ over $K$,  and a rational function $g \in K(C)^*$,
on associates to it its divisor $div_{C}(g) \in Z_{0}(C)$.
Given a  morphism $f: C \to X$, one may then consider the zero-cycle $f_{*}(div_{C}(g)) \in Z_{0}(X)$. 

One then defines the Chow group $CH_{0}(X)$ of zero-cycles on $X$ as the quotient of
$Z_{0}(X)$ by the subgroup spanned by all $f_{*}(div_{C}(g))$, for 
$f : C \to X$ a {\it proper} $K$-morphism from a normal, integral $K$-curve to $X$  and $g \in K(C)^*$.

 If $\phi : X \to Y$ is a proper morphism of $K$-varieties, there is an induced map
 $\phi_{*} : CH_{0}(X) \to CH_{0}(Y)$. In particular, if $X/K$ is proper, the structural map
 $X \to \Spec(K)$ induces a degree map $CH_{0}(X) \to \Z$.
 
 If $\phi : U \to X$ is an open embedding of $K$-varieties, the natural restriction map
 $Z_{0}(X) \to Z_{0}(U)$ which forgets closed points outside of $U$
 induces a map $CH_{0}(X) \to CH_{0}(U)$.

 Let $X$ be a $K$-variety. There is a natural bilinear pairing
 $$Z_{0}(X) \times \Br(X) \to \Br(K)$$
 which sends a pair $(P,\alpha)$ with $P$ a closed point of $X$ and
 $\alpha \in \Br(X)$ to $Cores_{K(P)/K}(\alpha(P))$.
  If $X/K$ is proper, 
  this pairing induces a bilinear pairing
 $$ CH_{0}(X) \times \Br(X) \to \Br(K).$$
  See \cite[Prop. 3.1]{ABBB}.

 This pairing satisfies an obvious functoriality property
 with respect to (proper) $K$-morphisms of proper $K$-varieties.

\section{Quadric surfaces over a field}\label{overafield}

The following proposition is classical. See \cite{Sk}, \cite{CTSD}, \cite[Thm. 3.1]{ACTP}.
 
\begin{prop} Let $K$ be a field, $char.(K) \neq 2$ and let $X \subset \P^3_{K}$ be a smooth quadric surface. It is defined by a quadratic form $q$, which one may assume to be  in diagonal
form $q= <1,-a,-b,abd>$, with $a,b,d \in K^*$. The class of  $d $ in $K^*/K^{*2}$ is the discriminant, it does not depend on the choice of the quadratic form $q$ defining the quadric $X$.

The natural map
$\Br(K) \to \Br(X)$ is surjective. 

(a) If $d \notin K^{*2}$, the map $\Br(K) \to \Br(X)$ is an isomorphism.

(b) If $d \in  K^{*2}$, the map $\Br(K) \to \Br(X)$ is surjective, and its kernel
is of order at most 2, spanned by the class of the quaternion algebra $(a,b)$,
which is nontrivial if and only if $X(K)=\emptyset$. 

\end{prop}

\section{A special quadric surface over $\P^2_{\C}  $}\label{specialexample}

Reference : \cite{HPT}, \cite{P}.

Let $F(x,y,z)= x^2+y^2+z^2-2(xy+yz+zx)$.

Let $X \subset \P^3_{\C} \times \P^2_{\C}$ be the family of 2-dimensional quadrics
over $\P^2_{\C}$
given by the  bihomogeneous equation
$$ yzU^2+zxV^2+xy W^2 + F(x,y,z) T^2=0.$$
This family is smooth over  the open set of $ \P^2_{\C}$ whose complement  
is the octic curve defined by the determinant equation
$$ \Delta= x^2.y^2.z^2. F(x,y,z)=0.$$
Note that this is the union of the smooth conic $F=0$ and
(twice) three tangents to this conic.
The family is flat over $\P^2_{\C}$ (all fibres are quadrics).
The total space is not smooth.

Part (a) of the following proposition is a result of Hassett, Pirutka, Tschinkel \cite[Prop. 11]{HPT}.

Part (b) is a special case of the general  statement \cite[Prop. 7]{Schr2},
the  proof of which  builds upon   results of Pirutka  (\cite[Thm. 3.17]{P}, \cite[Thm. 4]{Schr2}),
for which material is offered in Appendix  B below.

As we shall see, the proof  given for (a) in  \cite[Prop. 11]{HPT}  is easily modified to
simultaneously give a proof of (b).

The proposition suffices for the special case described in this note;
it dispenses us with the recourse to Appendix B.

\begin{prop}\label{calculcle}\label{HPT11}
Let $\tilde{X} \to X$ be a projective desingularisation of $X$.
Let $\alpha=(x/z,y/z) \in \Br(\C(\P^2))$.

(a)  The image $\beta$ of
$\alpha$  under the inverse map $p_{2}^* : \Br (\C(\P^2)) \to \Br(\C(X))$
is nonzero and lies in the subgroup $\Br(\tilde{X})$.

(b) For each codimension 1 subvariety $Y$ of $\tilde{X}$ which does not lie
over the generic point of $\P^2_{\C}$, the element $\beta \in \Br(\tilde{X})$
maps to $0 \in \Br(\C(Y))$.
\end{prop}

{\it Proof} 
The equation is symmetrical in $(x,y,z)$.
 The class $\alpha=(x/z,y/z)$ is given by $(x,y)$ in the open set $z\neq 0$,
by $(x/z,1/z)= (x,z)$ in the open set $y \neq 0$ and by $(1/z,y/z)=(y,z)$ in the open set $x\neq 0$.
In view of the symmetry between $(x,y,z)$ in the equation, we may restrict attention
to the open set $z\neq 0$. From now one we use affine coordinates $(x,y)$.
 In affine coordinates, the quaternion algebra $(x,y)$ has nontrivial residues
 along $x=0$ and $y=0$.

 Let $K=\C(\P^2)$. Let $X_{\eta}/K$ be the (smooth)  generic
quadric.
The discriminant of the quadratic form $q=<y,x,xy, F(x,y,1)> $ in $K^*$
is not a square. 
Thus the map $\Br(K)  \to \Br(X_{\eta})$ is an isomorphism (see \S \ref{overafield}).
Since
  the quaternion algebra $\alpha=(x,y) \in \Br(\C(\P^2))$ has  some nontrivial residues,
  it is nonzero in $\Br(\C(\P^2))$.
  Thus its image $\beta \in \Br(\C(X))$ is nonzero.

Let $v$ be a discrete valuation of rank one on $L:=K(X)$, let $S$
be its valuation ring. Let $\kappa_{v}$ denote the residue field.
 If $K \subset S$, then $(x,y)$ is unramified.
Suppose $S \cap K= R$ is a discrete valuation ring (of rank one).
The image of the closed point of $\Spec(R)$ in $\P^2_{\C}$
is then either a point  $m$ of codimension 1 or a (complex) closed point $m$ of $\P^2_{\C}$.

  Consider the first case. If the codimension 1 point $m$ does not belong
to $xy=0$, then $\alpha=(x,y) \in \Br(K) $ is unramified at $m$ on $\P^2_{\C}$
hence also in $\Br(L)$ at $v$. Moreover, the evaluation of $\beta$ in $\Br(\kappa_{v})$
 is just the image under $\Br(\C(m)) \to \Br(\kappa_{v})$ of the image of $\alpha$ in $\Br(\C(m))$,
 hence vanishes since
 $\Br(\C(m))=0$ (Tsen).

Suppose $m$ is a generic point of a component of $xyz=0$.
By symmetry,
 it is enough to examine the affine case
where the point $m$ of codimension 1 is the generic point of $x=0$.
In the function field $L$, we have an identity
$$yU^2+ x V^2  + xy W^2 + F(x,y,1)=0$$
with $yU^2+xV^2 \neq 0$.
In the completion of $K$ at the generic point of $x=0$,
$F(x,y,1)$ is a square, because $F(x,y,1)$ modulo $x$ is equal to $(y-1)^2$, a  nonzero square.
Thus in the completion $L_{v}$ we have an equality (with some other elements $U,V,W \in L_{v}$).
$$yU^2+ x V^2 +   xyW^2 + 1=0.$$ 
This gives $(x,y)_{L_{v}}=0 \in \Br(L_{v})$. Hence $(x,y)_{L}$ is unramified at $v$,
thus belongs to $\Br(R)$
and has image $0$ in $\Br(\kappa_{v})$.

  Suppose we are in the second case, i.e. $m$ is a closed point of $\P^2_{\C}$.
There is a  local map $O_{\P^2_{\C},m} \to S$ which induces a map $\C \to \kappa_{v}$.
  If $x\neq 0$, then
$x$ becomes a nonzero square in the residue field $\C$ hence in $\kappa_{v}$,
and the residue of $(x,y)_{L}$ at $v$ is trivial.
The analogous argument holds if $y\neq 0$. It remains to discuss the case
$x=y=0$.  We have $F(0,0,1)= 1 \in \C^*$. Thus $F(x,y,1)$ reduces to $1$
in $\kappa_{v} $, hence is a square in the completion $L_{v}$. 
As above, in the completion $L_{v}$
we have an equality
$$yU^2+ x V^2 +   xyW^2 + 1=0,$$ 
which implies $(x,y)_{L_{v}}=0 \in \Br(L_{v})$. 
Hence $(x,y)_{L}$ is unramified at $v$, thus belongs to $\Br(S)$
and has image $0$ in $\Br(\kappa_{v})$.
 $\Box$

\medskip

As in the reinterpretation \cite{CTO} of the Artin--Mumford examples,
 the intuitive idea behind the above result is that the quadric bundle  $X \to \P^2_{\C}$ is
ramified along $x.y.z.F(x,y,z)=0$ on $\P^2_{\C}$ and that 
the ramification of the symbol $(x/z,y/z)$ on $\P^2_{\C}$,
which is ``included'' in the ramification  of the quadric bundle $X \to \P^2_{\C}$
disappears over smooth projective models of $X$ : ramification eats up ramification
(Abhyankar's lemma). Here one also uses the fact that the smooth conic defined by $F(x,y,z)=0$
is tangent to each of the 
  lines $x=0, y=0, z=0$, and does not vanish at any of the intersection of
  these three lines.

\section{The specialisation argument}\label{specialisatonmethod}

 The following theorem is an improvement by S. Schreieder \cite[Prop. 26]{Schr1}
 of the specialisation method, as initiated by C. Voisin \cite{V}, in the format later
  proposed by  
 Colliot-Th\'el\`ene and Pirutka \cite{CTP}.  
 The assumptions in  \cite[Prop. 26]{Schr1} are more general than the ones
 given here. The generic fibre need not be smooth and one only requires that
 $f^{-1}(U) \to U$ be universally $CH_{0}$-trivial.
 There is a more general version which involves higher unramified cohomology
 with torsion coefficients. The proof is identical to the one given here  
 with the Brauer group.
 
 \medskip
 
  Schreieder's proof is cast in the geometric language of the decomposition of the diagonal.
 I provide a more ``field-theoretic'' proof. It is well known that both points of
 view are equivalent  \cite{ACTP,CTP}. I  add a further, hopefully simplifying, twist
 by using  specialization of $\rm R$-equivalence on rational points instead of Fulton's specialisation theorem
 for the Chow group.

\begin{theo}\label{main}
Let $R$ be a discrete valuation ring, $K$ its field of fractions, $\kappa$ its residue field.
Assume $\kappa$ is algebraically closed and $char(\kappa)=0$. Let $\overline{K}$ be
an algebraic closure of $K$.
Let $\mathcal{X}/R$ be an integral projective scheme over $R$,
with generic fibre $X=\mathcal{X}_{K}/K$ smooth, geometrically integral, and with special
fibre $Z/\kappa$ geometrically integral. 
Assume there exists a nonempty open set $U \subset Z$ and a projective, birational desingularisation
  $f : \tilde{Z} \to Z$ such that $V:=f^{-1}(U) \to U$ is an isomorphism, and such that
the complement $ \tilde{Z} \setminus V$ is a union $\cup_{i} Y_{i}$ of smooth irreducible divisors of $\tilde{Z}$.
Assume that the $\overline{K}$-variety $X_{\overline K}$
is  stably rational. If an element $\alpha \in \Br(\tilde{Z})$ vanishes on each $Y_{i}$,
then $\alpha=0 \in  \Br(\tilde{Z})$.
\end{theo}
{\it Proof}  
To prove the result, one may assume $R=\kappa[[t]]$ (completion of the original $R$)
and $K=\kappa((t))$.  Assume $X_{\overline K}$ is stably rational. Then there exists a finite extension $K_{1}=\kappa((t^{1/n}))$ of $K$
over which $X_{K_{1}}$ is $K_{1}$-stably rational. We replace $\mathcal{X}/R$ by
$\mathcal{X}\times_{R}   \kappa[[t^1/n]]$.   This does not change the special fibre.

Changing notation once more, we now have  $\mathcal{X}/R$ an integral projective scheme
whose generic fibre $X/K$  is  stably rational over $K$
and whose special fibre $Z/\kappa$ is just as in the theorem. 
Fix $m\in V(\kappa)$, mapping to $n \in U(\kappa)$.

Let $L=\kappa(Z)$.
We have the commutative diagram of exact sequences
$$\begin{array}{ccccccccc}
 &   &  \oplus_{i} CH_{0}(Y_{i,L}) & \to & CH_{0}(\tilde{Z}_{L}) & \to & CH_{0}(V_{L}) &\to & 0\\
&&  && \downarrow && \downarrow{\simeq} &&\\
  &  &   & &CH_{0}(Z_{L}) & \to & CH_{0}(U_{L}) &\to & 0.
 \end{array}$$
where  for each $i$, the closed embedding  $\rho_{i} = Y_{i} \to \tilde{Z}$ 
induces $$\rho_{i,*} :  CH_{0}(Y_{i,L}) \to CH_{0}(\tilde{Z}_{L}),$$ the top exact sequence is
the classical localisation sequence for the Chow group, the map 
$ f_{*}: CH_{0}(\tilde{Z}_{L}) \to CH_{0}(Z_{L})$ is induced by the proper map
$f: \tilde{Z} \to Z$, the map $ CH_{0}(V_{L}) \to CH_{0}(U_{L})$ is the isomorphism 
induced by the isomorphism\footnote{Instead of assuming that  $f^{-1}(U)  \to U$ is an isomorphism,
 it would be enough, as in \cite{Schr1},  to assume that this morphism 
 is a universal $CH_{0}$-isomorphism.} $f : V \to U$, and the map $CH_{0}(Z_{L}) \to CH_{0}(U_{L})$
is the obvious restriction map for the open set $U \subset Z$.

Let $\xi$ be the generic point of $\tilde{Z}$ and $\eta$ the generic point of $Z$.

Both $\eta_{L}$ and $n_{L}$ are smooth points of $Y_{L}$.
There exists an extension $R \to S$ of complete dvr inducing $\kappa \to L$
on residue fields.  Let $F$ be the field of fractions of $S$.
By Hensel's lemma, the points $\eta_{L}$ and $n_{L}$ lift to
rational points of the generic fibre  of $X\times_{K}F/F$ of $\mathcal{X}_{S}/S$.
Since $X/K$ is stably rational,  all points of $X_{F} (F)$ are 
${\rm R}$-equivalent (\cite[Prop. 10]{CTSansuc}, \cite[Cor. 6.6.6]{KS}).

It is a well known fact (\cite[prop. 3.1]{Madore}, \cite[Comments after Thm. 6.6.2]{KS}) that for a proper morphism
$\mathcal{X}_{S} \to S$ over a discrete valuation ring $S$ there is an induced map on ${\rm R}$-equivalence classes
$X(F)/{\rm R} \to Z(L)/{\rm R}$.  This implies $\eta_{L}-n_{L} = 0 \in CH_{0}(Z_{L})$.
\footnote{Alternatively, one could argue as follows.  
Since $X$ is stably rational over $K$, over any field $F$ containing $K$, the degree map $CH_{0}(X_{F}) \to \Z$
is an isomorphism (for a simple proof, see \cite[Lemme 1.5]{CTP}).  One could then invoke 
Fulton's specialisation theorem for the Chow group of a proper scheme over a dvr \cite[\S 2, Prop. 2.6]{F},
to get $\eta_{L}-n_{L} = 0 \in CH_{0}(Z_{L})$. Fulton's specialisation theorem is a nontrivial theorem.
The argument via ${\rm R}$-equivalence (cf. \cite[Remarque 1.19]{CTP}) looks  simpler.}

From the above diagram we conclude that
$$\xi_{L} = m_{L} + \sum_{i} {\rho_{i}}_{*}(z_{i}) \in CH_{0}(\tilde{Z}_{L})$$
with  $z_{i} \in  CH_{0}(Y_{i,L})$.

For the proper variety $\tilde{Z}_{L}$, there is a natural  bilinear pairing
$$ CH_{0}(\tilde{Z}_{L})  \times \Br(\tilde{Z}) \to \Br(L).$$
For the smooth, proper, integral  variety $\tilde{Z}$, on the generic point $\xi \in \tilde{Z}_{L}(L)$,
this pairing induces the embedding $\Br(\tilde{Z}) \hookrightarrow \Br(\kappa(Z))$.
Suppose $\alpha \in \Br(\tilde{Z})$    vanishes in each $\Br(Y_{i})$
(which follows from the vanishing in   $\Br(\kappa(Y_{i}))$ because $Y_{i}$ is smooth).
The evaluation of $\alpha$ on $m_{L}$ is just the image of $\alpha(m) \in \Br(\kappa)=0$.
The above equality implies $\alpha(\xi)= 0 \in \Br(L)$, hence $\alpha=0 \in \Br(\tilde{Z})$.
$\Box$.

\section{Stable rationality is not constant in smooth projective families}\label{stabnotstab}

 We now complete the simplified proof of the theorem of Hassett, Pirutka and Tschinkel \cite{HPT}.
 
 \medskip

\begin{theo}\label{mainmain}
There exist a smooth projective family of complex 4-folds $f : X \to T$
parametrized by an open set $T$ of the affine line $\A^1_{\C}$ and points
$m, n \in T(\C)$ such that the fibre $X_{n}$ is rational  and the fibre $X_{m}$
is not stably rational.
\end{theo}

{\it Proof}
One  considers the universal family of quadric bundles over $ \P^2_{C}$
given
 in $ \P^3_{\C} \times \P^2_{C}$
  by a bihomogenous form of bidegree $(2,2)$. This is given by a symmetric $(4,4)$ square matrix  with entries $a_{i,j}(x,y,z)$ homogeneous quadratic forms in three variables $(x,y,z)$.
If its determinant is nonzero, it is
  a homogeneous polynomial of degree 8.

We thus have a parameter space $B$  given 
 by a projective space of dimension 59 (the corresponding vector space
 being  given by  the coefficients of 10 quadratic forms in three variables).
 We have the map $X\to B$ whose fibres $X_{m}$ are the various  
 quadric bundles $X_{m} \to \P^2_{\C}$, for $X_{m}\subset \P^3_{\C} \times \P^2_{C}$
given by the vanishing of a  nonzero complex bihomogeneous form of bidegree $(2,2)$.

Using Bertini's theorem, one shows that
there exists a nonempty open set $B_{0} \subset B$ such that
the fibres of $X\to B$ over points of $m\in B_{0}$ are
 flat quadric bundles   $X_{m} \to \P^2_{\C}$
 which are smooth as $\C$-varieties.

 Using Bertini's theorem, one also shows that there exist
 points $m\in B_{0}$ with the property that the corresponding
quadric bundle has 
  $a_{1,1} =0$, which implies that the fibration $X_{m} \to \P^2_{\C}$
  has a rational section (given by the point $(1,0,0,0)$), hence that
  the generic fibre of $X_{m} \to \P^2_{\C}$ is rational over $\C(\P^2)$,
  hence that the $\C$-variety $X_{m}$ is rational over $\C$.
  [Warning : this  Bertini argument uses the fact that we consider families
 of quadric surfaces over $\P^2_{\C}$. It does not work for families
 of conics over $\P^2_{\C}$.]
 
 These Bertini arguments are briefly described in \cite[Lemma 20 and  Thm. 47]{Schr1} and are tacitly used   in
 \cite[Page 3]{Schr2}.
 
By Proposition \ref{calculcle},
the special example in \S  \ref{specialexample} defines a point  $P_{0} \in B(\C)$ 
 whose fibre is $Z=X_{P_{0}}$ and which admits a  projective birational desingularisation
 $f : \tilde{Z} \to Z$ satisfying :
 
(a) there exists a nonempty open set $U \subset Z$, 
 such that the induced map  $V:=f^{-1}(U) \to U$ is an isomorphism;
 
 (b) the complement $ \tilde{Z} \setminus V$ is a union $\cup_{i} Y_{i}$ of smooth irreducible divisors of $\tilde{Z}$;

(b) there is a nontrivial element $\alpha \in \Br(\tilde{Z})$  which vanishes on each $Y_{i}$.

Theorem \ref{main} then implies that the generic fibre of $X \to B$ is not geometrically 
stably rational.  There are various ways to conclude from this that there are many points
$m\in B_{0}(\C)$ such that the fibre $X_{m}$ is not stably rational.

Take one such point $m\in B_{0}(\C)$ and a point  $n\in B_{0}(\C)$ such that
$X_{n}$ is rational. Over an open set of the line joining $m$ and $n$
we get a projective family of  smooth varieties with one fibre rational
and with one fibre not stably rational.
$\Box$

\medskip

 The proof by Hassett, Pirutka and Tschinkel \cite{HPT}  uses an explicit
desingularisation of the variety $Z$ in \S \ref{specialexample},
with a description of the exceptional divisors  appearing in the process. 
Schreieder's improvement of the specialisation method enables one to bypass this explicit desingularisation.

Note that   \cite{HPT} and \cite{Schr2} contain many more
results on families of quadrics surfaces over $\P^2$ than  Theorem \ref{mainmain}.

\section{Appendix A. Conics over  a discrete valuation ring}

Let $R$ be a dvr with residue field $k$ of characteristic not 2. Let $K$ be the fraction field.
A smooth conic over $K$ admits a regular model $\mathcal{X}$ given in $\P^2_{R}$ either by an equation
$$x^2-ay^2-bz^2=0$$
with $a,b \in R^*$ (case (I))
or  a regular model  $\mathcal{X}$ given by an equation
$$x^2-ay^2-\pi z^2=0$$
with $a\in R^*$ and $\pi$ a uniformizing parameter (case (II)).
Moreover, in the second case one may assume that $a$
is not a square in the residue field $\kappa$.

\begin{prop}\label{conicoverDVR}
Let $R$ be a dvr with residue field $k$ of characteristic not 2. Let $K$ be the fraction field.
Let $W \to \Spec(R)$ be a proper flat morphism with $W$ regular and connected.
Assume that the generic fibre is a smooth conic over $K$.
Then 

(a) The natural map $\Br(R) \to \Br(W)$ is onto.

(b) For $Y \subset W$  an integral divisor
 contained in the special fibre of $W \to \Spec(R)$,
 and $\beta \in \Br(W)$,
 the  image of $\beta$ under restriction $\Br(W) \to \Br(Y)$ belongs to the image
 of $\Br(\kappa)  \to \Br(Y)$.
\end{prop}

{\it Proof}
By purity for the Brauer group of a 2-dimensional regular scheme, to prove (a),
one may assume that $W=\mathcal{X}$ as above. Let $X=\mathcal{X}\times_{R}K$.
It is well known that the map $\Br(K) \to \Br(X)$ is onto,
with kernel spanned by the quaternion symbol $(a,b)_{K}$ in case (I)
and by $(a,\pi)_{K}$ in case (II).

Let $\beta\in \Br(\mathcal{X}) \subset \Br(X)$. 
Let $\alpha \in \Br(K)$ be some element with image $\beta_{K}$.
We have the exact sequence 
$$0 \to \Br(R)\{2\} \to \Br(K)\{2\} \to H^1(\kappa,\Q_{2}/\Z_{2})$$
Comparison of residues on $\Spec(R)$ and on $\mathcal{X}$
 shows that  the residue $\delta_{R}(\alpha)$ is either 0
 or is equal to the nontrivial class in $H^1(k(\sqrt{\overline a})/k,\Z/2)$,
 and this last case may happen only in case (II).
 In the first case, we have $\alpha\in \Br(R)$, hence $\beta-\alpha_{\mathcal{X}}=0$
 in $\Br(X_{K})$ hence also in $\Br(X)$ since $X$ is regular. In the second case, we have
  $$\delta_{R}(\alpha)=\delta_{R}((a,\pi))$$
  hence $\alpha = (a,\pi) + \gamma$ with $\gamma \in \Br(R)$.
  We then get $$\beta= (a,\pi)_{K(X)} + \gamma_{K(X)} \in \Br(K(X)).$$
  But $(a,\pi)_{K(X)}=0$. Thus $\beta-\gamma_{\mathcal{X}} \in  \Br(\mathcal{X}) \subset
  \Br(K(X))$ vanishes, hence $\beta= \gamma_{\mathcal{X}}  \in \Br(\mathcal{X})$.
  The map $\Br(R) \to \Br(\mathcal{X})$ is thus surjective. This gives (a) for $(\mathcal{X}$
  hence for $W$,
  and (b) immediately follows.
  $\Box$

\bigskip

{\it Exercise} Artin-Mumford type examples are specific singular conic bundles $X$  in the total space of
a rank 3 projective bundle over $\P^2_{\C}$
whose unramified Brauer group is non trivial. 
Using Proposition \ref{conicoverDVR} and Theorem \ref{main}, deform 
 such examples into conic bundles of the same type with smooth ramification locus
 and whose total space is not stably rational. As in section \ref{specialexample}, there is no need to compute
 an explicit resolution of singularities of $X$.
 l

\section{Appendix B. Quadric surfaces  over a discrete valuation ring}\label{overadvr}

The following section was written up to give details on some tools
and results used in \cite[Thm. 4]{Schr2}. As demonstrated above, this section turns out not to be necessary
to vindicate the HPT example. But it is useful   for more general examples.

References : \cite{Sk}, \cite[\S 3]{CTSko}, 
  \cite[Thm. 2.3.1]{CTSD},  \cite[Thm. 3.17]{P}.

\medskip

Let $R$ be a discrete valuation ring, $K$ its fraction field,
$\pi$ a uniformizer, $\kappa=R/(\pi)$ the residue field.
Assume $char(\kappa) \neq 2$.

Let $X \subset \P^3_{K}$ be a smooth quadric, defined by a nondegenerate
4-dimensional quadratic form $q$.  Up to scaling and changing of variables,
there are four possibilities.

(I) $q=<1,-a,-b,abd>$ with $a,b,d \in R^*$.

(II) $q=<1,-a,-b,\pi>$ with $a, b \in R^*$ and $\pi$ a uniformizing parameter of $R$.

(III) $q=<1,-a, \pi, -\pi.b>$ with $a, b  \in R^*$ and $\pi$ a uniformizing parameter of $R$.
The class of $a.b \in R^*$ represents the discriminant of the quadratic form.
Its image $\overline{a}.\overline{b} \in \kappa^*$ is a square if and  only if the discriminant
of $q$ is a square in the completion of $K$ for the valuation defined by $R$.

Let ${\mathcal X} \subset \P^3_{R}$ be the subscheme
cut out by $q$. Let $Y/\kappa$ be the special fibre.

In case (I),  ${\mathcal X} /R$ is smooth. 

In case (II), $X$ is regular, 
the special fibre $Y$ is a cone over a smooth conic.

In case (III), the special fibre is given by the equation
$x^2-\overline{a}y^2=0$ in $\P^3_{\kappa}$.
If $\overline{a}$ is a square, this is the union of two planes
intersection along the line $x=y=0$. If If $\overline{a}$ is not a square,
this is an integral scheme which over $\kappa(\sqrt{a})$ breaks up
as the union of two planes. In both cases, the scheme $\mathcal{X}$
is singular at the points $x=y=0, z^2-\overline{d}t^2=0$.
See \cite[\S 2]{Sk}.

\begin{prop}
Let us assume $char(\kappa)=0$.

In case (III), let $W \to \mathcal{X}$ be a projective, birational desingularisation of $\mathcal X$.

In case (I), the map $\Br(R) \to \Br(\mathcal{X})$ is onto.
If $d \in R$ is not a square, it is an isomorphism. If $d$ 
is a square, the kernel is spanned by the class $(a,b) \in \Br(R)$.

In case (II), the map $\Br(R) \to \Br(\mathcal{X})$ is an isomorphism.

In case (III), assume   $\overline{a}.\overline{b}$ is not a square in $\kappa$.
Then  $\Br(R) \to \Br(W)$ is onto.

In case (III), 
if either $\overline{a}$ or $\overline{b}$ is a square, or if
 $\overline{a}.\overline{b}$ is not a square, then $\Br(R) \to \Br(W)$ is onto.
 An element of $\Br(K)$ whose image in $\Br(X)$ lies in $ \Br(W)$
 belongs to $\Br(R)$.
 
 In case (III),   assume   $\overline{a}.\overline{b}$
is a square in $\kappa$.  Then the image of $(a,\pi) \in \Br(K)$ in $\Br(X)$
 belongs to $\Br(W)$. It spans the quotient of $\Br(W)$ by the image of
 $\Br(R)$.
 If  moreover $\overline{a}$ is not a square in $\kappa$,
then  it does not belong to the image of $\Br(R)$. 
 \end{prop}

 {\it Proof}
 Let $x$ be 
 a codimension 1 regular point on  $\mathcal{X}$ or on $W$,
 lying above the closed point of $\Spec(R)$. Let $e_{x}$ denote its multiplicity in the fibre.
 We have a commutative diagram
 $$\begin{array}{ccccccccc}
  \Br(K) & \to & H^1(\kappa,\Q/\Z)\\
  \downarrow & & \downarrow \\
  \Br(X) &\to &H^1(\kappa(x), \Q/\Z)\\
  \end{array}
  $$
  The kernel of $Br(K) \to H^1(\kappa,\Q/\Z)$ is $\Br(R)$.
  
  In case (I) and (III), the special fibre $Y$ is geometrically integral over $\kappa$,
  the multiplicity is $1$, the map $H^1(\kappa,\Q/\Z) \to H^1(\kappa(x), \Q/\Z)$
  is thus injective.  This is enough to prove the claim.
  
 Let us consider case (III). The map $\Br(K) \to \Br(X)$ is onto. Let $\alpha \in \Br(K)$.
 Let $\rho \in H^1(\kappa,\Q/\Z)$ be its residue.
 On the (singular) normal model given by   $q=<1,-a, \pi, -\pi.b>$ over $R$, 
 if $\overline{a} \in \kappa$ is a square, the fibre $Y$ contains   geometrically integral components
 of multiplicity $1$  given by the components of $x^2-\overline{a}y^2=0$.
 By the above diagram,  $\rho=0 \in H^1(\kappa,\Q/\Z)$.
 We can also use the model given by  $q=<1,-b, \pi, -\pi.a>$. If $\overline{b} \in \kappa$ is a square, we conclude that $\rho=0 \in H^1(\kappa,\Q/\Z)$. Let us assume that $\rho\neq 0 \in 
 H^1(\kappa,\Q/\Z)$. Thus $\overline{a} $ and $\overline{b}$ are nonsquares. On the first model,
 the kernel of $H^1(\kappa, \Q/\Z) \to H^1(\kappa(Y),\Q/\Z)$ coincides with the kernel
 of $H^1(\kappa, \Q/\Z) \to H^1(\kappa(\sqrt{\overline{a}}, \Q/\Z)$, which is the
 $\Z/2$-module spanned by the class of $\overline{a}$ in $\kappa^*/\kappa^{*2}=H^1(\kappa,\Z/2)$. On the second model, the kernel of $H^1(\kappa, \Q/\Z) \to H^1(\kappa(Y),\Q/\Z)$ 
 is the
 $\Z/2$-module spanned by the class of $\overline{b}$ in $\kappa^*/\kappa^{*2}=H^1(\kappa,\Z/2)$.  We thus conclude that $\overline{a}.\overline{b}$  is a square in $\kappa$, and that
 the residue of $\alpha$ coincides with $\overline{a}$, i.e. is equal to the residue of
 $(a,\pi) \in \Br(K)$ (or to the residue of $(b,\pi)$).
 
It remains to show that if  $\overline{a}.\overline{b}$ is a square in $\kappa$, then
 $(a,\pi)$ has trivial residues on $W$ and more generally
with respect to any rank one discrete valuation $v$ on the function field $K(X)$ of $X$.
One may restrict attention to those $v$  which induce the $R$-valuation on $K$.
Let $S \subset K(X)$ be the valuation ring of $v$ and let $\lambda$ be its residue field.
 There is an inclusion $\kappa \subset \lambda$. In $K(X)$ we have  an equality
 $$(x^2-ay^2)=\pi.(z^2-b),$$
 where both sides are nonzero. Thus
 in $\Br(K(X))$, we have the equality
 $$(a,\pi)= (a,x^2-ay^2) + (a, z^2-b) = (a,z^2-b),$$
 where the last equality comes from the classical $(a,x^2-ay^2)=0$.
 To compute residues, we may go over to completions.
 In the completion of $R$, $a.b$ is a square. It is thus a square in the completion of
 $K(X)$ at $v$. But then in this completion $(a,z^2-b)=(b,z^2-b)=0$
 Hence the residue of $(a,\pi)$ at $v$ is zero.
 $\Box$
 
\begin{prop}
Assume $char(\kappa)=0$.
Let ${\mathcal{X}} \subset \P^3_{R}$ be as above, and let $W \to \mathcal{X}$
be a proper birational map with $W$ regular.
Let $\beta \in \Br(W)$ and let $Y \subset W$ be an integral divisor
 contained in the special fibre of $W \to \Spec(R)$.
Then the image of $\beta$ in $\Br(\kappa(Y))$ belongs to the image
of $\Br(\kappa) \to \Br(\kappa(Y))$.
\end{prop}

 {\it Proof}
In case (I) and (II), and in case (III) when $\overline{a}.\overline{b}$ is a square in $\kappa$,
this is clear since then the map $\Br(R) \to \Br(W)$ is onto.

Suppose we are in case (III). To prove the result, we may make a base change from $R$
to its henselisation.
 Then $ab$ is square in $R$. The group $\Br(W)$ is spanned  by the image of $\Br(R)$ and the image of
the class $(a,\pi)$. The equation of the quadric may now be written
$$ X^2-aY^2+\pi Z^2- a\pi T^2=0.$$ 
This implies that $(a,-\pi)$ vanishes in the Brauer group of the function field $\kappa(W)$ of $W$. Since $W$ is regular,
the map $\Br(W) \to \Br(\kappa(W))$ is injective. Since $(a,-\pi)_{\kappa(W)}$ belongs to $\Br(W)$ and
spans $\Br(W)$ modulo the image of $\Br(R)$, this completes the proof. 
$\Box$

\medskip

One may rephrase the above results in a simpler fashion.
\begin{prop}
Assume $char(\kappa)=0$.
Let ${\mathcal{X}} \subset \P^3_{R}$ be as above, and let $W \to \mathcal{X}$
be a proper birational map with $W$ regular.

(i) If $R$ is henselian, then the map 
$\Br(R) \to \Br(W)$ is onto.

(ii) For any element $\beta \in \Br(W)$ and
$Y \subset W$   an integral divisor
 contained in the special fibre of $W \to \Spec(R)$,
 the  image of $\beta$ under restriction $\Br(W) \to \Br(Y)$ belongs to the image
 of $\Br(\kappa)  \to \Br(Y)$.
\end{prop}

Upon use of Merkurjev's geometric lemmas \cite[\S 1]{M}, and use of Tsen's theorem,
one then gets
\cite[Prop. 7]{Schr2} of Schreieder.

\section{Appendix C. A remark on the vanishing of unramified elements on components of the special fibre}

The following proposition, found in June 2017,
 gives some partial explanation for the vanishing on components
of the special fibre which occurs in 
  \cite[Prop. 6, Prop. 7]{Schr1} \cite[Prop. 7]{Schr2} or in Proposition \ref{calculcle} above.
Unfortunately the proof requires that the component be of multiplicity one in the fibre.
Since this was written, in the case of quadric bundles, S. Schreieder \cite[\S 8.3]{Schr3}
has managed to use arguments as in  \cite[\S 3]{CTSko}
to get information on what happens with the other components.

\medskip

\begin{prop} Let  $A \hookrightarrow B$ be a local homomorphism of discrete valuation rings
and let 
 $K \subset L$ be the inclusion of their fraction fields.
Let  $\kappa \subset \lambda$ be the induced inclusion on their residue fields.

Let $\ell$ be a prime invertible in $A$.

Let $i \geq 2$ be an integer and let  $\alpha \in H^{i}(K,\mu_{\ell}^{\otimes i})$. 

Assume:

(i) $B$ is unramified over $A$.

(ii) The image of $\alpha$ in $H^{i}(L,\mu_{\ell}^{\otimes i})$ is unramified,
and in particular is the image of a (well defined) element 
 $\beta \in H^{i}(B,\mu_{\ell}^{\otimes i})$.

Then $\beta(\lambda) \in H^{i}(\lambda,\mu_{\ell}^{\otimes i})$
is in the image of  $H^{i}(\kappa, \mu_{\ell}^{\otimes i}) \to H^{i}(\lambda,\mu_{\ell}^{\otimes i})$.
\end{prop}

{\it Proof}   I  use Gersten's conjecture for $H^{i}$ with finite coefficients
over discrete valuation rings. I also use Bloch-Kato's conjecture as proven by
Rost, Voevodsky, Weibel. Going over to henselisations would probably dispense with 
the appeal to such deep results.

Let  $\pi \in A$ be a uniformizer.
By the Bloch-Kato  conjecture, the class
 $\alpha \in H^{i}(K,\mu_{\ell}^{\otimes i})$ is a sum of symbols.
The standard relation $(a,-a)=0$ for $a \in K^*$ implies that $\alpha$
 is the sum of     $\alpha_{1} \in H^{i}(A,\mu_{\ell}^{\otimes i})
\subset H^{i}(K,\mu_{\ell}^{\otimes i})$
and  
$\alpha_{2}=(-\pi,a_{2}, \dots, a_{i})$ with all $a_{j} \in A^*$.

By hypothesis, the image of  $(-\pi,a_{2}, \dots, a_{i})$ dans $ H^{i}(L,\mu_{\ell}^{\otimes i})$ is
unramified, hence  it is equal to some $\beta_{2} \in H^{i}(B,\mu_{\ell}^{\otimes i})$.

Since $B$ is unramified over $A$, the uniformizer  $\pi$ is also a uniformizer of  $B$.
To compute the image  $\beta_{2}(\lambda)$  of $\beta_{2} \in H^{i}(B,\mu_{\ell}^{\otimes i})$
in  $H^{i}(\lambda,\mu_{\ell}^{\otimes i})$, one may use the formula
$$\beta_{2}(\lambda)= \partial_{B}(\pi, \beta_{2}) \in H^{i}(\lambda,\mu_{\ell}^{\otimes i}).$$

Now
  $$\beta_{2}(\lambda)= \partial_{B}(\pi, -\pi,a_{2}, \dots, a_{i})=0$$
since $(\pi,-\pi)=0$. 
Thus  $\beta(\lambda) \in  H^{i}(\lambda,\mu_{\ell}^{\otimes i})$
coincides with the image of  
$\alpha_{1}(\kappa) \in H^{i}(\kappa,\mu_{\ell}^{\otimes i})$
under the map
$H^{i}(\kappa, \mu_{\ell}^{\otimes i}) \to H^{i}(\lambda,\mu_{\ell}^{\otimes i})$.
$\Box$

\end{document}